\documentclass[12pt]{article}
\usepackage{amsfonts}
\usepackage{latexsym,amssymb,amsfonts}
\usepackage{mathrsfs}
\usepackage{graphicx}
\usepackage{psfrag}
\usepackage{amsmath, amsbsy}
\usepackage{amsopn, amstext}
\usepackage{amsmath,multicol,amsopn}
\usepackage{latexsym,amssymb}

\allowdisplaybreaks[3]

\def\sqr#1#2{{\vcenter{\vbox{\hrule height.#2pt
              \hbox{\vrule width.#2pt height#1pt \kern#1pt \vrule width.#2pt}
              \hrule height.#2pt}}}}
\def\signed #1{{\unskip\nobreak\hfil\penalty50
              \hskip2em\hbox{}\nobreak\hfil#1
              \parfillskip=0pt \finalhyphendemerits=0 \par}}
\def\endpf{\signed {$\sqr69$}}

\def\dbR{{\mathop{\rm l\negthinspace R}}}

\def\3n{\negthinspace \negthinspace \negthinspace }
\def\2n{\negthinspace \negthinspace }
\def\1n{\negthinspace }

\def\dbR{{\mathbb{R}}}

\def\ds{\displaystyle}

\def\={\buildrel \triangle \over =}

\def\resp{{\it resp. }}
\def\a{\alpha}

\def\d{\delta}
\def\e{\varepsilon}
\def\z{\zeta}

\def\l{\lambda}

 \def\n{\nabla}

\def\t{\times}
\def\f{\varphi}
\def\th{\theta}
\def\o{\omega}

\def\ns{\noalign{\ss} }

\def\diag{\hbox{\rm $\,$diag$\,$}}

\def\pa{\partial}
\def\G{\Gamma}

\def\Si{\Sigma}

\def\O{\Omega}

\def\mE{\mathbb{E}}

\def\cA{{\cal A}}

\def\cE{{\cal E}}
\def\cF{{\cal F}}

\def\cM{{\cal M}}

\def\cO{{\cal O}}

\def\no{\noindent}

\def\ms{\medskip}
\def\bs{\bigskip}
\def\q{\quad}
\def\qq{\qquad}

\def\max{\mathop{\rm max}}
\def\min{\mathop{\rm min}}

\def\pa{\partial}

\def\wt{\widetilde}
\def\cd{\cdot}
\def\cds{\cdots}

\def\div{\hbox{\rm div$\,$}}

\def\dist{\hbox{\rm dist$\,$}}

\def\|{\Big |}
\def\({\Big (}
\def\){\Big )}
\def\[{\Big[}
\def\]{\Big]}
\def\be{\begin{equation}}
\def\bel{\begin{equation}\label}
\def\ee{\end{equation}}
\def\bt{\begin{theorem}}
\def\bcd{\begin{condition}}
\def\ecd{\end{condition}}
\def\et{\end{theorem}}
\def\bc{\begin{corollary}}
\def\ec{\end{corollary}}
\def\bde{\begin{definition}}
\def\ede{\end{definition}}
\def\bl{\begin{lemma}}
\def\el{\end{lemma}}
\def\bp{\begin{proposition}}
\def\ep{\end{proposition}}
\def\br{\begin{remark}}
\def\er{\end{remark}}
\def\ba{\begin{array}}
\def\ea{\end{array}}
\def\ed{\end{document}}
\def\ns{\noalign{\ms}}
\def\ds{\displaystyle}

\def\square#1{\vbox{\hrule\hbox{\vrule height#1%
     \kern#1\vrule}\hrule}}
\def\rectangle#1#2{\vbox{\hrule\hbox{\vrule height#1%
     \kern#2\vrule}\hrule}}

\font\tenbb=msbm10 \font\sevenbb=msbm7 \font\fivebb=msbm5

\newfam\bbfam
\scriptscriptfont\bbfam=\fivebb \textfont\bbfam=\tenbb
\scriptfont\bbfam=\sevenbb

\newtheorem{lemma}{Lemma}[section]
\newtheorem{remark}{Remark}[section]

\newtheorem{theorem}{Theorem}[section]
\newtheorem{corollary}{Corollary}[section]

\newtheorem{definition}{Definition}[section]
\newtheorem{proposition}{Proposition}[section]
\newtheorem{condition}{Condition}[section]

\newtheorem{lm}{Lemma}[section]

\makeatletter
   
   \@addtoreset{equation}{section}
\makeatother

\begin{document}
\title{\bf  Observability Estimate and  State Observation Problems  for Stochastic Hyperbolic Equations\thanks{This work is partially supported by the NSF of China under grants
 11101070. This paper is an improved
version of one chapter of the author's Ph D thesis (\cite{Luqi2})
accomplished at Sichuan University under the guidance of Professor
Xu Zhang. The author would like to take this opportunity to thank
him deeply for his help. \ms}}

\author{Qi L\"u\thanks{Universit\'{e} Pierre et Marie Curie-Paris VI
UMR 7598, Laboratoire Jacques-Louis Lions, 4, place Jussieu
Paris, F-75005 France; School of Mathematical Sciences, University of Electronic
 Science and Technology of China, Chengdu, 610054, China.
  {\small\it E-mail:} {\small\tt luqi59@163.com}.}  }

\date{}

\maketitle

\begin{abstract}\no
In this paper, we derive a boundary and
an internal observability inequality
for stochastic hyperbolic equations
with nonsmooth lower order terms. The
required inequalities are obtained by
global Carleman estimate for stochastic
hyperbolic equations. By these
inequalities, we study a state
observation problem for stochastic
hyperbolic equations. As a consequence,
we also establish a unique continuation
property for stochastic hyperbolic
equations.
\end{abstract}

\bs

\no{\bf 2010 Mathematics Subject Classification}.  Primary  65N21, 60H15, 93B07.

\bs

\no{\bf Key Words}. Stochastic hyperbolic
equation, observability, state observation problem, Carleman estimate,
unique continuation property.

\ms

\section{Introduction }
\q Let $T > 0$, $G \in \mathbb{R}^{n}$ ($n
\in \mathbb{N}$) be a given bounded domain
with the $C^{2}$ boundary $\G$. Let  $\G_0$
be a suitable chosen nonempty subset of
$\G$, whose definition will be given later.
Put $$
\begin{array}{ll}\ds
Q \= (0,T) \t G,\q \Si \= (0,T) \t \G
,\q \Si_0 \= (0,T) \t \G_0,\\
\ns\ds \cO_\d(\G_0) \= \Big\{ x\in G :\,
\dist(x,\G_0)\leq \d \Big\} \mbox{ for some }\d>0.
\end{array}
$$

Let $(\O, {\cal F}, \{{\cal F}_t\}_{t \geq 0},
P)$ be a complete filtered probability space on
which a  one dimensional standard Brownian motion
$\{ B(t) \}_{t\geq 0}$ is defined. Let $H$ be a
Banach space. Denote by $L^{2}_{\cal F}(0,T;H)$
the Banach space consisting of all $H$-valued and $\{
{\cal F}_t \}_{t\geq 0}$-adapted processes
$X(\cdot)$ such that
$\mathbb{E}(|X(\cdot)|^2_{L^2(0,T;H)}) < \infty$,
by $L^{\infty}_{\cal F}(0,T;H)$ the Banach space
consisting of all $H$-valued and $\{ {\cal F}_t
\}_{t\geq 0}$-adapted bounded processes, by
$L^{2}_{\cal F}(\O;C([0,T];H))$ the Banach space
consisting of all $H$-valued and $\{ {\cal F}_t
\}_{t\geq 0}$-adapted processes $X(\cdot)$ such
that $\mathbb{E}(|X(\cdot)|^2_{C(0,T;H)}) <
\infty$(similarly, one can define $L^{2}_{\cal
F}(\O;C^{k}([0,T];H))$ for any positive integer
$k$), all of these spaces are endowed with the
canonical norm .

Throughout this paper, we make the following
assumptions on the coefficients  $ b^{ij} \in
C^1(G)$:

1.\;\,$b^{ij} = b^{ji}$ $(i,j = 1,2,\cdots, n)$;

2. For some constant $s_0
> 0$,
\begin{equation}\label{bij}
 \sum_{i,j}b^{ij}\xi^{i}\xi^{j} \geq s_0 |\xi|^2,
\,\,\,\,\,\,\,\,\,\, \forall\, (x,\xi)\=
(x,\xi^{1}, \cdots, \xi^{n}) \in G \t
\mathbb{R}^{n}.
\end{equation}
 Here and in what follows, we denote $\ds \sum_{i,j = 1}^{n}$ simply by $\ds
\sum_{i,j}$. For simplicity, we  use the notation
$\ds y_i \equiv y_{i}(x) \= \frac{\partial
y(x)}{\partial x_i}$, where $x_i$ is the $i$-th
coordinate of a generic point $x=(x_1,\cdots,
x_n)$ in $\mathbb{R}^{n}$. In a similar manner,
we  use notations $z_i$, $v_i$, etc. for the
partial derivatives of $z$ and $v$ with respect
to $x_i$. Also, we denote by $\nu(x) = (\nu^1(x),
\cdots, \nu^n(x))$ the unit outward normal vector of
$\G$ at point $x$.

\vspace{0.1cm}

Let us consider the following stochastic
hyperbolic equation:
\begin{eqnarray}{\label{system1}}
\left\{
\begin{array}{lll}\ds
\ds dz_{t} -
\sum_{i,j}(b^{ij}z_i)_{j}dt = \big[b_1
z_t +
 b_2\cd\nabla z  + b_3 z + f \big]dt + (b_4 z + g)dB(t) & {\mbox {
in }} Q,
 \\
\ns\ds  z = 0 & \mbox{ on } \Si, \\
\ns\ds  z(0) = z_0, z_{t}(0) = z_1 & \mbox{ in }
G.
\end{array}
\right.
\end{eqnarray}
Here the initial data $(z_0, z_1) \in L^2(\O,{\cal
F}_0, P; H^1_0(G) \t L^2(G))$, the coefficients $b_i$ $(1 \leq i \leq 4)$
satisfy  that
\begin{equation} \label{aibi}
\begin{array}{ll}\ds
  b_1 \in L_{\cal
F}^{\infty}(0,T;L^{\infty}(G)), \qq  b_2 \in
L_{\cal
F}^{\infty}(0,T;L^{\infty}(G;\mathbb{R}^{n})),  \\
\ns\ds b_3 \in L_{\cal
F}^{\infty}(0,T;L^{p}(G))\,(p\in [n,\infty]), \qq
\,\, b_4 \in L_{\cal
F}^{\infty}(0,T;L^{\infty}(G)),
\end{array}
\end{equation}
and nonhomogeneous terms
\begin{eqnarray}\label{fg}
\qq f \in L^2_{\cF}(0,T;L^2(G)), \qq g \in
L^2_{\cF}(0,T;L^2(G)).
\end{eqnarray}

Put  \begin{equation}\label{HT}
H_{T} \= L_{\cal F}^2 (\O; C([0,T];H_{0}^1(G)))\cap
L_{\cal F}^2 (\O; C^{1}([0,T];L^2(G))).
\end{equation}
Clearly, $H_T$ is a Banach space with the
canonical norm.

\vspace{0.1cm}

Now we give the  definition
of  the solution to the equation
 \eqref{system1}.

\begin{definition} \label{def solution to sys}
We call $z\in H_T$  a solution to the equation
\eqref{system1} if the following two conditions hold: \\1.
$z(0) = z_0$ in $G$, P-a.s., and
$z_t(0) = z_1$ in $G$, P-a.s.  \\
2. For any $t \in [0,T]$ and any $\eta \in H_0^1(G)$, it holds
that
\begin{equation} \label{solution to sysh}
\begin{array}{ll}\ds
 \q \int_{G} z_t(t,x)\eta(x)dx - \int_{G} z_t(0,x)\eta(x)dx
  \\ \ns\ds= \int_0^t \int_G \Big\{
-\sum_{i,j}b^{ij}(x)z_i(s,x)\eta_j(x) + \big[b_1(s,x) z_t (s,x) +
 b_2(s,x)\cd\nabla z(s,x)\\
 \ns\ds \qq + b_3(s,x) z(s,x)  + f(s,x)\big]\eta(x) \Big\}dxds   \\
\ns\ds \q + \int_0^t \int_G \big[ b_4(s,x)
z(s,x) + g(s,x)\big]\eta(x) dxdB(s), \,\, \mbox { P-a.s.}
\end{array}
\end{equation}
\end{definition}

For any initial data $(z_0, z_1) \in L^2(\O,{\cal
F}_0, P; H^1_0(G) \t L^2(G))$, one can show that the
equation \eqref{system1} admits a unique solution
$z \in H_T $(see \cite{Zhangxu3} for details).

\vspace{0.2cm}

Before giving $\G_0$, we introduce the following
condition:
\begin{condition}
\label{condition of d}
There exists a positive
function $d(\cdot) \in
C^2(\overline{G})$ satisfying the following:\\
{\rm 1}. For some constant $\mu_0 > 0$, it holds
\begin{equation}\label{d1}
\begin{array}{ll}\ds
\sum_{i,j}\Big\{ \sum_{i',j'}\Big[
2b^{ij'}(b^{i'j}d_{i'})_{j'} -
b^{ij}_{j'}b^{i'j'}d_{i'} \Big]
\Big\}\xi^{i}\xi^{j} \geq \mu_0
\sum_{i,j}b^{ij}\xi^{i}\xi^{j}, \\
\ns\ds \hspace{5.5cm} \forall\,
(x,\xi^{1},\cdots,\xi^{n}) \in  \overline{G}  \t
\mathbb{R}^n.
\end{array}
\end{equation}
{\rm 2}. There is no critical point of $d(\cdot)$
in $\overline{G}$, i.e., \be\label{d2}\min_{x\in
\overline{G} }|\nabla d(x)| > 0. \ee
\end{condition}

\begin{remark}
If $(b^{ij})_{1\leq i,j\leq n}$ is the identity
matrix, then $d(x)=|x-x_0|^2$ satisfies Condition
\ref{condition of d}, where $x_0$ is any point
which belongs to $\mathbb{R}^n\setminus\overline
G$.
\end{remark}

\begin{remark}\label{rm1}
Condition \ref{condition of d} was first given in
\cite{Fu-Yong-Zhang1} for the purpose of
obtaining an internal observability estimate for
hyperbolic equations. In that paper, the authors
also gave some explanation of Condition
\ref{condition of d} and some interesting
nontrivial examples satisfying it. Further, a detailed study of this condition is given in \cite{Liu}.
\end{remark}

The $\G_0$ is as follows:
\begin{eqnarray}\label{def gamma0}
\G_0 \= \Big\{ x\in \G \Big| \sum_{i,j}b^{ij}d_i(x)\nu^{j}(x) > 0
\Big\}.
\end{eqnarray}

It is easy to check that if $d(\cdot) $
satisfies Condition \ref{condition of
d}, then for any given constants $a
\geq 1$ and $b \in \mathbb{R}$, the
function $\tilde{d} = ad + b$ still
satisfies Condition  \ref{condition of
d}  with $\mu_0$ replaced by $a\mu_0$.
Therefore  we may choose $d$, $\mu_0$,
$c_0>0$, $c_1>0$ and $T$ to such
that  the following  condition holds:
\begin{condition}\label{condition2}
\begin{equation}\label{con2 eq1}
1. \qq  \frac{1}{4}\sum_{i,j}b^{ij}(x)d_i(x)d_j(x) \geq R^2_1\=\max_{x\in\overline G}d(x)\geq R_0^2 \= \min_{x\in\overline G}d(x),\q \forall x\in \overline G.
\end{equation}
\qq\; {\rm 2}. \q $T> T_0\=2 R_1.$\\

\q\,\!\! {\rm 3}. \q $\ds\(\frac{2R_1}{T}\)^2<c_1<\frac{2R_1}{T}$.\\

\q\,\!\!  {\rm 4}. \q $\mu_0 - 4c_1 -c_0 > 0$.
\end{condition}

\begin{remark}\label{rm2}
As we have explained, since
$\ds\sum_{i,j}b^{ij}d_id_j >0$, and one
can choose $\mu_0$ in Condition
\ref{condition of d} large enough,
Condition \ref{condition2}   could be
satisfied obviously. We put it here
just in order to emphasize the
relationship between $0< c_0< c_1 <1$,
$\mu_0$ and $T$.
\end{remark}
\begin{remark}\label{rmT}
If $(b^{ij})_{1\leq i,j\leq n}$ is the
identity matrix, then it is easy to
show that $$d(x)=2|x-x_0|^2$$ for some
$x_0\notin \overline G$ satisfy
\eqref{d1} and \eqref{d2} in Condition
\ref{condition of d}. However, this $d(\cd)$
does not satisfy \eqref{con2 eq1} in
the Condition \ref{condition2}. On the
other hand, if we consider the problem
with $(b^{ij})_{1\leq i,j\leq
n}=\diag(1,1,\cds,1)$, we do not need
\eqref{con2 eq1}. Indeed, in
this case, the inequality \eqref{obser
esti2} and \eqref{inobser esti2} below hold
for all $\ds T> 2\max_{x\in\overline
G}|x-x_0|$. One can follow the proofs
of Theorem \ref{observability} and
\ref{inobser} to see this. We omit the
details.
\end{remark}

In the rest of this paper, we use $C$
to denote a generic positive constant
depending on $G$, $T$, $\G_0$,
$b^{ij}$($i,j=1,\cds,n$), $d$, $c_0$
and $c_1$(unless otherwise stated),
which may change from line to line.

Put
\begin{eqnarray}\label{r1r2}
r_1\= |b_2|_{L_{\cal
F}^{\infty}(0,T;L^{\infty}(G;\mathbb{R}^{n}))}
+ |(b_1,b_4)|_{L_{\cal
F}^{\infty}(0,T;(L^{\infty}(G))^2}\;\mbox{ and
}\; r_2 =   |b_3|_{L_{\cal
F}^{\infty}(0,T;L^{p}(G))}.
\end{eqnarray}

Now we give our main results.
The first one is the boundary observability
estimate for the equation \eqref{system1}.

\begin{theorem}\label{observability}
Let Condition  \ref{condition of d} and Condition
\ref{condition2} be satisfied. For any solution
of the equation
 \eqref{system1}, we have
\begin{equation} \label{obser esti2}
\begin{array}{ll}\ds
\q |(z_0,z_1)|_{L^2(\O,{\cal F}_0, P; H_0^1(G)\t
L^2(G))}
\\ \ns\ds  \leq Ce^{C(r_1^2 + r_2^{\frac{1}{ 3/2 - n/p}}+1)} \Big(\Big|\frac{\partial z}{\partial \nu}\Big |_{L^2_{\cal
 F}(0,T;L^2(\G_0))} + |f|_{L^2_{\cal
 F}(0,T;L^2(G))} + |g|_{L^2_{\cal
 F}(0,T;L^2(G))}\Big).
\end{array}
\end{equation}
\end{theorem}
The second one is   the internal
observability estimate for the equation
\eqref{system1}.
\begin{theorem}\label{inobser}
Let Condition \ref{condition of d} and
Condition \ref{condition2}  be
satisfied. For any solution of the equation
 \eqref{system1}, it holds
\begin{equation} \label{inobser esti2}
\begin{array}{ll}\ds
\q|(z_0,z_1)|_{L^2(\O,{\cal F}_0, P; H_0^1(G)\t
L^2(G))}
\\ \ns\ds \leq  e^{C(r_1^2 + r_2^{\frac{1}{ 3/2 - n/p}}+1)} \Big(\big|\nabla z\big |_{L^2_{\cal
 F}(0,T;L^2(\cO_\d(\G_0)))} + |f|_{L^2_{\cal
 F}(0,T;L^2(G))} + |g|_{L^2_{\cal
 F}(0,T;L^2(G))}\Big).
\end{array}
\end{equation}
\end{theorem}

\begin{remark}
Inequality  \eqref{obser esti2}
{\resp(\eqref{inobser esti2})}  is
referred to as  observability estimate
since it provides a quantitative
estimate of the norm of the initial
data in terms of the observed quantity,
by means of the observability constant
$C$. Indeed, the inequality \eqref{obser
esti2} {\resp(\eqref{inobser esti2})}
allows one to estimate the total energy
of solutions at time $0$ in terms of
the partial energy localized in the
observation subboundary
$\Gamma_0$(\resp the observation
subdomain $\cO_{\d}$). This sort of
inequality is strongly relevant to
control problems and state observation problems for stochastic
hyperbolic equations.
\end{remark}

\begin{remark}\label{rmI}
Compared with inequality \eqref{inobser
esti2}, it is  more interesting  to
establish  the following inequality:
\begin{equation} \label{inobser esti2.1}
\begin{array}{ll}\ds
\q|(z_0,z_1)|_{L^2(\O,{\cal F}_0, P; L^2(G)\t
H^{-1}(G))}
\\ \ns\ds \leq e^{C(r_1^2 + r_2^{\frac{1}{ 3/2 - n/p}}+1)} \Big(| z  |_{L^2_{\cal
 F}(0,T;L^2(\cO_\d(\G_0)))} + |f|_{L^2_{\cal
 F}(0,T;L^2(G))} + |g|_{L^2_{\cal
 F}(0,T;L^2(G))}\Big).
\end{array}
\end{equation}
However, we do not know how to obtain this result now.
\end{remark}

Thanks to its important applications in
Control Theory and Inverse Problems for
hyperbolic equations, and to its strong
connection with the unique continuation
for solutions to hyperbolic equations,
the observability estimate for
hyperbolic equations have been studied
extensively in the literature. There
are four main approaches to study it.
The first one is the multiplier
techniques(see \cite{Lions1} for
example). The second one is nonharmonic
Fourier series technique(see
\cite{Komornik1} for example). The
third one is based on the Microlocal
Analysis(see
\cite{Bardos-Lebeau-Rauch1} for
example). The last one is the global
Carleman estimate(see
\cite{DZZ,Fu-Yong-Zhang1} for example).

Among the above four methods, the global
Carleman estimate is the most common
and powerful technique to derive
observability inequalities. It can be
regarded as a more developed version of
the classical multiplier technique.
Compared with  the multiplier method,
the Carleman approach is robust with
respect to the lower order terms.
Compared with the microlocal analysis,  it
requires less regularity on
coefficients and domain. Compared with
the nonharmonic Fourier series method,
it has much less restrictions on the
shape of the domain.

There are very few works addressing the
observability problems for stochastic
partial differential equations. To the
best of our knowledge,
\cite{barbu1,Luqi4,Tang-Zhang1,Zhangxu3}
are the only references for this topic.
In \cite{barbu1,Tang-Zhang1} the
observability estimate for  stochastic
heat equations is studied. \cite{Luqi4}
is devoted to the observability
estimate for stochastic Schr\"{o}dinger
equations while \cite{Zhangxu3} is
concerned with the observability
estimate for stochastic wave equations.
In \cite{Zhangxu3}, the author proves a
boundary observability estimate for the
equation \eqref{system1} with
$(b^{ij})_{1\leq i,j \leq n}$ being the
identity matrix. More precisely, the
author proves that
\begin{equation} \label{obser esti22}
\begin{array}{ll}\ds
\q|(z(t),z_t(t))|_{L^2(\O,{\cal F}_t,
P; H_0^1(G)\t L^2(G))}
\\ \ns\ds \leq Ce^{Ct^{-1}\cA} \Big(\Big|\frac{\partial z}{\partial \nu}\Big |_{L^2_{\cal
 F}(0,T;L^2(\G_0))} + |f|_{L^2_{\cal
 F}(0,T;L^2(G))} + |g|_{L^2_{\cal
 F}(0,T;L^2(G))}\Big),
\end{array}
\end{equation}
where $z$ solves the equation
\eqref{system1}, $T$ satisfies
\begin{equation}\label{T}
\frac{\ds(4+5c)\min_{x\in\overline G}|x-x_0|^2}{9c} > c^2 T^2 > 4\max_{x\in\overline G}|x-x_0|^2
\end{equation}
for some $c\in (0,1)$ and $x_0\in
\mathbb{R}^n\setminus\overline G$, and
$$
\cA = \cA(b_1,b_2,b_3,b_4)\=
|b_2|^2_{L_{\cal
F}^{\infty}(0,T;L^{\infty}(G;\mathbb{R}^{n}))}
+ |(b_1,b_4)|^2_{L_{\cal
F}^{\infty}(0,T; L^{\infty}(G))^2} +
|b_3|^2_{L_{\cal
F}^{\infty}(0,T;L^{n}(G))}.
$$

There are three main differences
between the inequality \eqref{obser
esti2} and \eqref{obser esti22}. The
first is that the left-hand side of
\eqref{obser esti2} is
$|(z_0,z_1)|_{H_0^1(G)\times L^2(G)}$.
However, it seems that one cannot
simply replace
$|(z(t),z_t(t))|_{L^2(\O,\cF_t,P;H_0^1(G)\times
L^2(G))}$ by
$|(z_0,z_1)|_{L^2(\O,\cF_0,P;H_0^1(G)\times
L^2(G))}$ directly in \eqref{obser
esti22}, thanks to the term
$e^{Ct^{-1}\cA}$. Although by
Proposition \ref{energy ensi} in this
paper, one can get the estimate for
$|(z_0,z_1)|_{L^2(\O,\cF_0,P;H_0^1(G)\times
L^2(G))}$ by \eqref{obser
esti22}, there are other two differences. The second is that the
observation time $T$ in \eqref{obser
esti22} should satisfy \eqref{T}, which
is usually much more restrictive than
that $T>2\max_{x\in\overline
G}|x-x_0|$ for our result(see
Remark~\ref{rmT}). The third is that
the observability constant  in
\eqref{obser esti22} is not as sharp as
that in \eqref{obser esti2}. Indeed, it
is clear that
$$
|b_1|_{L^{\infty}(Q)}^2+|
b_2|_{L^\infty(Q;\,\mathbb{R}^{n})}^2+|b_3|_{L^{\infty}(0,T;\,L^p(G))}^{\frac{1}{3/2-n/p}}\le
\cA(b_1,b_2,b_3,0)+C\qq\forall\;p\in[n,\infty].
$$

The rest of this paper is organized as
follows. In Section 2, we collect some
preliminaries. Section 3 is devoted to the proof of Theorem \ref{observability}. In Section 4, we prove Theorem \ref{inobser}. Section 5 is addressed to a state observation problem of semilinear stochastic hyperbolic equations. At last, in Section 6, we present some further comments and open problems.

\section{Some Preliminaries }

In this section, we present some preliminary
results. First, we give a hidden regularity property of solutions to the
equation \eqref{system1}.
\begin{proposition}\label{hidden r}
For any solution of the equation \eqref{system1}, it
holds that
\begin{equation}\label{hidden ine}
\begin{array}{ll}\ds
\q\Big|\frac{\partial z}{\partial \nu}\Big
|_{L^2_{\cal
 F}(0,T;L^2(\G_0))} \\
 \ns\ds \leq e^{C(r_1^2 + r_2^2 + 1)} \Big(|(z_0,z_1)|_{L^2(\O,{\cal F}_0, P;
H_0^1(G)\t L^2(G))} + |f|_{L^2_{\cal
 F}(0,T;L^2(G))} + |g|_{L^2_{\cal
 F}(0,T;L^2(G))}\Big).
 \end{array}
\end{equation}
\end{proposition}
\begin{remark}
In \cite{Zhangxu3}, the author proved Proposition
\ref{hidden r} for the case $(b^{ij})_{1\leq i,\leq
j\leq n}$ is identity matrix. The proof
Lemma \ref{hidden r} for general
$(b^{ij})_{1\leq i, j\leq n}$ is similar. We
only give a sketch of it here.
\end{remark}

{\it Proof}\,:  For any $$h\=(h^1,\cds,h^n)\in
C^1(\dbR_t \t \dbR^n_x ; \dbR^n), $$ by direct
computation, we can show that
\begin{equation}\label{equality hidden1}
\begin{array}{ll}
\ds \q-\sum_{i=1}^n\Big[ 2(h\cd\nabla z)\sum_{j=1}^n b^{ij}z_{x_j} +
h^i\Big( z_t^2 - \sum_{i,j=1}^n
b^{ij}z_{x_i} z_{x_j} \Big) \Big]_{x_i}dt\\
\ns =  \ds 2 \Big[\Big(dz_t - \sum_{i,j=1}^n
(b^{ij}z_{x_i})_{x_j}dt\Big) h \cd \nabla z -
d(z_t h\cd \nabla z) + z_t h_t\cd \n zdt -
\sum_{i,j,k=1}^n b^{ij}z_{x_i} z_{x_k}
h^k_{x_j}dt\Big] \\
\ns \ds\q  - (\div h) z_t^2dt + \sum_{i,j=1}^n z_{x_j} z_{x_i}
\div(b^{ij}h)dt.
\end{array}
\end{equation}
Since $\G\in C^2$, one can find a vector field
$\xi=(\xi^1,\cdots,\xi^n)\in
C^1(\mathbb{R}^n,\mathbb{R}^n)$ such that $\xi=\nu$
on $\G$(see \cite[page 18]{Komornik}). Setting
$h=\xi$ in the equality \eqref{equality hidden1},
integrating  \eqref{equality hidden1} in
$Q$, taking expectation in $\O$ and  integrating
by parts, we get inequality \eqref{hidden ine}
immediately. \endpf

\vspace{0.2cm}

Further, we give an energy estimate for the equation
\eqref{system1}, which plays an important role in
the proof of the  observability estimate.

\medskip

\begin{proposition}\label{energy ensi}
For any $z$ solves the equation \eqref{system1}, it
holds that
\begin{equation}\label{en esti}
\begin{array}{ll}
\q\ds \mE\int_G\big( |z_t(t,x)|^2 + |\nabla
z(t,x)|^2 \big) dx \\
\ns \leq \ds  Ce^{C\(r_1^2 + r_2^{\frac{1}{2-n/p}}+1\)T} \mE\int_G
\big(|z_t(s,x)|^2 + |\nabla z(s,x)|
\big) dx
\\
\ns \ds \q+\;
 C\mE\int_0^T\int_G\big[ f^2(\tau,x) + g^2(\tau,x)
\big]dxd\tau \Big\},
 \end{array}
\end{equation}
for any $0\leq s,t\leq T$.
\end{proposition}

\medskip

{\it Proof}\,:  Without loss of generality, we
assume that $t \leq s$. Let
$$
\cE(t) = \mE\int_G\[ |z_t(t,x)|^2 + |\nabla
z(t,x)|^2 +  r_2^{\frac{2}{2-n/p}} |z(t,x)|^2\]
dx.
$$
From Poincar\'{e}'s inequality, we get
\begin{equation}\label{en eq0}
\mE\int_G\big( |z_t(t)|^2 + |\nabla z(t)|^2 \big)
dx \leq \cE(t) \leq C\(r_2^{\frac{2}{2-n/p}}
+1\)\mE\int_G\big( |z_t(t)|^2 + |\nabla z(t)|^2
\big) dx.
\end{equation}
By means of It\^{o}'s formula, we have
$$
d(z^2_t) = 2z_t dz_t + (dz_t)^2,
$$
which implies that
\begin{equation}\label{en eq1}
\begin{array}{ll}
\ds \q\mE\int_G  \(|z_t(s,x)|^2 + r_2^{\frac{2}{2-n/p}}|z(s,x)|^2\) dx -  \mE\int_G \(|z_t(t,x)|^2 + r_2^{\frac{2}{2-n/p}}|z(t,x)|^2\)   dx  \\
\ns\ds =  -2\mE\int_t^s\int_G \sum_{i,j}b^{ij}(x)(\tau,x)z_i(\tau,x) z_{jt}(\tau,x)dxd\tau  + \mE\int_t^s\int_G z_t(\tau,x)\Big[b_1(\tau,x) z_t(\tau,x) \\
\ns \ds \q +
2 b_2(\tau,x)\cd\nabla z(\tau,x)  + b_3(\tau,x) z(\tau,x) + f(\tau,x) \Big]dxd\tau \\
\ns \ds \q + \mE\int_t^s\int_G \big[b_4(\tau,x) z(\tau,x)
 + g(\tau,x) \big]^2 dxd\tau + 2r_2^{\frac{2}{2-n/p}}\mE\int_t^s\int_G z_t(\tau,x)z(\tau,x)dxd\tau.
\end{array}
\end{equation}
Therefore, we obtain that
\begin{equation}\label{en eq2}
\begin{array}{ll}
\ds \q\mE\int_G  \[|z_t(s,x)|^2 + \sum_{i,j}b^{ij}(x)z_i(s,x)z_j(s,x) \] dx \\
\ns \ds \q  -  \mE\int_G  \[|z_t(t,x)|^2 + \sum_{i,j}b^{ij}(x)z_i(t,x)z_j(t,x) \] dx\\
\ns \ds = \mE\int_t^s\int_G z_t(\tau,x)\Big[b_1(\tau,x)
z_t(\tau,x)+
 b_2(\tau,x)\cd\nabla z(\tau,x)  + b_3(\tau,x) z(\tau,x) + f (\tau,x)\Big]dxd\tau \\
\ns  \ds \q + \mE\int_t^s\int_G \big[b_4(\tau,x) z(\tau,x)
 + g(\tau,x) \big]^2 dxd\tau+ 2r_2^{\frac{2}{2-n/p}}\mE\int_t^s\int_G z_t(\tau,x)z(\tau,x)dxd\tau\\
\ns\ds \leq C(r_1^2 \!+\! 1)\mE\!\int_t^s\!\!\int_G\!
\big[z^2_t(\tau,x)\! +\! |\nabla z(\tau,x)|^2 \!+\!
z^2(\tau,x)\big]dxd\tau \!+\! 2r_2^{\frac{2}{2-n/p}}\mE\!\!\int_t^s\!\int_G z_t(\tau,x)z(\tau,x)dxd\tau \\
\ns \ds \q  + \mE\int_t^s\int_G b_3(\tau,x)
z(\tau,x)z_t(\tau,x)dxd\tau + 2\mE\int_t^s\int_G\big[ f^2(\tau,x) +
g^2(\tau,x) \big]dxd\tau.
\end{array}
\end{equation}

Put $p_1 = \frac{2p}{n-2}$ and $p_2 =
\frac{2p}{p-n}$. It is easy to check that
$$
\frac{1}{p} + \frac{1}{p_1} + \frac{1}{p_2} +
\frac{1}{2} = 1\;\;\mbox{ and }\;
\frac{1}{2(n/p)^{-1}} + \frac{1}{2(1-n/p)^{-1}} +
\frac{1}{2} = 1.
$$
By H\"{o}lder's inequality and Sobolov's
embedding theorem, we find
\begin{equation}\label{en eq2.1}
\begin{array}{ll}
\ds \q\Big|\mE\int_{G}b_3(\tau,x) z(\tau,x)z_t(\tau,x)dx\Big|
\\
\ns\ds \leq\mE\int_{G}|b_3(\tau,x)|
|z(\tau,x)|^{\frac{n}{p}}|z(\tau,x)|^{1-\frac{n}{p}}|z_t(\tau,x)|dx\\
\ns\ds \leq r_2\mE
\(\big||z(\tau,\cd)|^{\frac{n}{p}}\big|_{L^{p_1}(G)}
\big||z(\tau,\cd)|^{1-\frac{n}{p}}\big|_{L^{p_2}(G)}\big|
z_t(\tau,\cd) \big|_{L^2(G)} \)\\
\ns \ds = r_2 \mE \( \big|
z(\tau,\cd)\big|^{\frac{n}{p}}_{L^{\frac{n}{n-2}}(G)}
 \big|
z(\tau,\cd)\big|^{1-\frac{n}{p}}_{L^{2}(G)}\big|z_t(\tau,\cd)
\big|_{L^2(G)} \)\\
\ns \ds = r_2^{\frac{1}{2-n/p}}\mE \( \big|
z(\tau,\cd)\big|^{\frac{n}{p}}_{L^{\frac{n}{n-2}}(G)}
r_2^{\frac{1-n/p}{2-n/p}}\big|
z(\tau,\cd)\big|^{1-\frac{n}{p}}_{L^{2}(G)}\big|z_t(\tau,\cd)
\big|_{L^2(G)} \).
\end{array}
\end{equation}
Since
$$
\left\{
\begin{array}{ll}\ds
\big|
z(\tau,\cd)\big|^{\frac{n}{p}}_{L^{\frac{n}{n-2}}(G)}
\leq \[\int_G\( |z_t(\tau,x)|^2 + |\nabla z(\tau,x)|^2
+ r_2^{\frac{2}{2-n/p}} |z(\tau,x)|^2\) dx
\]^{\frac{n}{2p}},\\
\ns\ds r_2^{\frac{1-n/p}{2-n/p}}\big|
z(\tau,\cd)\big|^{1-\frac{n}{p}}_{L^{2}(G)} \leq
\[\int_G\( |z_t(\tau,x)|^2 + |\nabla z(\tau,x)|^2 +
r_2^{\frac{2}{2-n/p}} |z(\tau,x)|^2\) dx
\]^{\frac{1}{2}-\frac{n}{2p}},\\
\ns\ds \big|z_t(\tau,\cd) \big|_{L^2(G)} \leq
\[\int_G\( |z_t(\tau,x)|^2 + |\nabla z(\tau,x)|^2 +
r_2^{\frac{2}{2-n/p}} |z(\tau,x)|^2\) dx
\]^{\frac{1}{2}},
\end{array}
\right.
$$
from \eqref{en eq2.1}, we get
\begin{equation}\label{en eq2.2}
\Big|\mE\int_{G}b_3(\tau,x)
z(\tau,x)z_t(\tau,x)dx\Big| \leq
r_2^{\frac{1}{2-n/p}}\cE(\tau).
\end{equation}
By a similar argument, we obtain that
\begin{equation}\label{en eq2.3}
\begin{array}{ll}\ds
\q r_2^{\frac{2}{2-n/p}}\mE\int_G
z(\tau,x)z_t(\tau,x)dx \\
\ns\ds \leq
\frac{1}{2}r_2^{\frac{1}{2-n/p}}\mE\int_G\(
r_2^{\frac{2}{2-n/p}}z^2(\tau,x) +
z_t^2(\tau,x)\)dx \\
\ns\ds \leq
\frac{1}{2}r_2^{\frac{1}{2-n/p}}\cE(\tau).
\end{array}
\end{equation}

From \eqref{en eq2.2}, \eqref{en eq2.3} and the
property of $b^{ij}$($i,j=1,\cds,n$)(see
\eqref{bij}), we find that
\begin{equation}\label{en eq3}
\begin{array}{ll}
\ds  \cE(t)   \leq C \Big\{\cE(s) +  \(r_1^2 +
r_2^{\frac{1}{2-n/p}}+1\) \int_t^s \cE(\tau)d\tau
+
 \mE\int_t^s\int_G\big[ f^2(\tau,x) + g^2(\tau,x)
\big]dxdt\Big\}.
\end{array}
\end{equation}
This, together with backward Gronwall's inequality,
implies that
\begin{equation}\label{en eq4}
\cE(t) \leq e^{C\(r_1^2 +
r_2^{\frac{1}{2-n/p}}+1\)(s-t)}\cE(s) +
C\mE\int_t^s\int_G\big[ f^2(\tau,x) + g^2(\tau,x)
\big]dxdt.
\end{equation}
From \eqref{en eq0} and \eqref{en eq4}, we get
\begin{equation}
\begin{array}{ll}
\ds \q\mE\int_G  \big(|z_t(t,x)|^2 + |\nabla z(t,x)|^2
\big) dx \\
\ns\ds \leq Ce^{C\(r_1^2 +
r_2^{\frac{1}{2-n/p}}+1\)(s-t)} \mE\int_G
\big(|z_t(s,x)|^2 + |\nabla z(s,x)|^2
\big) dx
\\
\ns\ds \qq +
 C\mE\int_t^s\int_G\big[ f^2(\tau,x) + g^2(\tau,x)
\big]dxdt \Big\},
\end{array}
\end{equation}
which leads to the inequality \eqref{en esti}
immediately.
\endpf

\vspace{0.2cm}

At last,  we introduce the following known result, which plays a key
role in getting the boundary and internal Carleman estimate.

\begin{lm}\cite[Theorem 4.1]{Zhangxu3}\label{hyperbolic1}
Let $p^{ij} \in C^{1}((0,T)\t \mathbb{R}^n)$ satisfy
\begin{equation}
p^{ij} = p^{ji}, \qq i,j = 1,2,\cdots,n,\nonumber
\end{equation}
$l,\Psi \in C^2((0,T)\t\mathbb{R}^n)$. Assume
that $u$ is an $H^2_{loc}(\mathbb{R}^n)$-valued
and $\{\cF_t\}_{t\geq 0}$-adapted process such
that $u_t$ is an $L^2(\mathbb{R}^n)$-valued
semimartingale. Set $\theta = e^l$ and $v=\theta
u$. Then, for a.e. $x\in \mathbb{R}^n$ and P-a.s.
$\o \in \O$,
\begin{eqnarray}\label{hyperbolic2}
&\,&\q\theta \Big( -2l_t v_t +
2\sum_{i,j}p^{ij}l_i v_j + \Psi v \Big)
\Big[ du_t - \sum_{i,j}(p^{ij}u_i)_j dt \Big] \nonumber\\
&\,& \,\,\,\q+\sum_{i,j}\Big[
\sum_{i'j'}\big(
2p^{ij}p^{i'j'}l_{i'}v_iv_{j'} -
p^{ij}p^{i'j'}l_iv_{i'}v_{j'}
\big) - 2p^{ij}l_t v_i v_t + p^{ij}l_i v_t^2 \nonumber\\
&\,& \qq \qq \qq \,\,\,\,+ \Psi p^{ij}v_i v - \Big( Al_i +
\frac{\Psi_i}{2}\Big)p^{ij}v^2 \Big]_j dt \\
&\,& \,\,\,\q +d\Big[ \sum_{i,j}p^{ij}l_t
v_i v_j - 2\sum_{i,j}p^{ij}l_iv_jv_t + l_t
v_t^2 - \Psi v_t v + \Big( Al_t +
\frac{\Psi_t}{2}\Big)v^2 \Big] \nonumber \\
&\,& = \Big\{ \Big[ l_{tt} + \sum_{i,j}(p^{ij}l_i)_{j} - \Psi
\Big]v_t^2 - 2\sum_{i,j}[(p^{ij}l_j)_t +
p^{ij}l_{tj}]v_iv_t \nonumber \\
&\,& \,\,\,\,\,\,\,\,+\sum_{i,j} \Big[ (p^{ij}l_t)_t +
\sum_{i',j'}\Big(2p^{ij'}(p^{i'j}l_{i'})_{j'}-(p^{ij}p^{i'j'}l_{i'})_{j'}\Big)
+ \Psi p^{ij} \Big]v_iv_j \nonumber \\
&\,&\,\,\,\,\,\,\,\,+ Bv^2 + \Big( -2l_tv_t +
2\sum_{i,j}p^{ij}l_iv_j + \Psi v \Big)^2\Big\} dt +
\theta^2l_t(du_t)^2, \nonumber
\end{eqnarray}
where $(du_t)^2$ denotes the quadratic
variation process of $u_t$, $A$ and $B$ are
stated as follows:
\begin{eqnarray}\label{AB1}
\left\{
\begin{array}{lll}\ds
A\=(l_t^2 - l_{tt}) - \sum_{i,j}(p^{ij}l_il_j -p_j^{ij}l_i -
p^{ij}l_{ij})-\Psi, \\
\ns\ds B\=A\Psi + (Al_t)_t-\sum_{i,j}(Ap^{ij}l_i)_j +
\frac{1}{2}\Big[ \Psi_{tt} - \sum_{i,j}(p^{ij}\Psi_i)_j \Big].
\end{array}
\right.
\end{eqnarray}
\end{lm}


\section{Proof of Theorem \ref{observability}}


In this section, we prove Theorem
\ref{observability} by means of the global
Carleman estimate. In what follows, for $\l \in
\mathbb{R}$, we  use $O(\l^r)$ to denote a
function of order $\l^r$ for large $\l$.

\vspace{0.2cm}

{\it Proof of Theorem \ref{observability}}\,: We
divide the proof into three steps.

\vspace{0.2cm}

{\bf Step 1}.
 On one hand,
from Condition \ref{condition2}, we know that
there is an $\e_1\in (0,1/2)$ such that
\begin{equation}\label{e1}
l(t,x)\leq \l\(\frac{R_1^2}{2} - \frac{cT^2}{8}\)
< 0,\q\forall\, (t,x)\in \[ \(0, \frac{T}{2}-\e_1
T\)\bigcup \(\frac{T}{2}-\e_1 T,T\)  \]\times G.
\end{equation}
On the other hand, since
$$
l\( \frac{T}{2},x\) = d(x)\geq R_0^2, \qq \forall x\in G,
$$
we can find an $\e_0\in (0,\e_1)$ such that
\begin{equation}\label{e0}
l(t,x)\geq \frac{R_0^2}{2},\q\forall (t,x)\in  \(\frac{T}{2}-\e_0 T,\frac{T}{2}+\e_0 T\) \times G.
\end{equation}

Now we choose a $\chi\in C^\infty_0[0,T]$ satisfying
\begin{equation}\label{chi}
\chi=1 \mbox{ in } \(\frac{T}{2}-\e_1 T,\frac{T}{2}+\e_1 T\).
\end{equation}
 Let $y=\chi z$ for $z$ solving the equation \eqref{system1}, then we know that $y$ is a solution to the following equation:
\begin{eqnarray}{\label{system2}}
\left\{
\begin{array}{lll}\ds
\ds dy_{t} - \sum_{i,j}(b^{ij}y_i)_{j}dt = \Big[b_1 y_t +
(b_2,\nabla y) + b_3 y + \chi f + \a \Big]dt + (b_4 y + \chi g)dB(t) & {\mbox {
in }} Q,
 \\
\ns\ds  y = 0 & \mbox{ on } \Si, \\
\ns\ds  y(0) = y(T) = 0, y_{t}(0) = y_t(T) = 0 & \mbox{ in } G.
\end{array}
\right.
\end{eqnarray}
Here $\a = \chi_{tt}z + 2\chi_t z_t - b_1\chi_t
z$.

\vspace{0.2cm}

{\bf Step 2.} We apply Lemma \ref{hyperbolic1} to the solution of the
equation  \eqref{system2}. In the present case,
we choose
$$
p^{ij} = b^{ij}, \Psi = l_{tt} +
\sum_{i,j}(b^{ij}l_i)_{j} - \l c_0,
$$
and then
estimate the terms in \eqref{hyperbolic2} one by
one.

 We first analyze  the terms which
stand for the ``energy" of the solution. The point is to compute the
order of $\l$ in the coefficients of $|v_t|^2$,
$|\nabla v|^2$ and $|v|^2$. Clearly, the term for
$|v_t|^2$ reads
\begin{equation}\label{coeffvt} \Big\{l_{tt} +
\sum_{i,j}(b^{ij}l_i)_{j} -\Psi\Big\}v_{t}^2 = \l c_0 v_{t}^2.
\end{equation}
Noting that $b^{ij}$($1 \leq i, j \leq n$)
are independent of $t$ and
$l_{tj}=l_{jt}=0$, we get that
\begin{equation}\label{bcoeffvtvi}
\sum_{i,j}\big[(b^{ij}l_j)_t +
b^{ij}l_{tj}\big]v_i v_t = 0.
\end{equation}
By Condition \ref{condition of d} and Condition
\ref{condition2}, we have that
\begin{equation}\label{vivj}
\begin{array}{ll}\ds
\q\sum_{i,j}\Big\{ (b^{ij}l_t)_t + \sum_{i',j'}\big[
2b^{ij'}(b^{i'j}l_{i'})_{j'} - (b^{ij}b^{i'j'}l_{i'})_{j'} \big] +
\Psi b^{ij} \Big\}v_i v_j \nonumber\\
\ns\ds = \sum_{i,j} \Big\{ 2b^{ij}l_{tt} - b^{ij}\l c_0 +
\sum_{i',j'}\big[ 2b^{ij'}(b^{i'j}l_{i'})_{j'} -
b^{ij}_{j'}b^{i'j'}l_{i'} \big] \Big\}v_i v_j \\
\ns\ds \geq \l (\mu_0 -4c_1 -
c_0)\sum_{i,j}b^{ij}v_i v_j.
\end{array}
\end{equation}
Now we compute the coefficients of $|v|^2$.
\begin{equation}
\begin{array}{ll}\ds
 A \3n &\ds= l_t^2 - l_{tt} - \sum_{i,j} \big[ b^{ij}l_i l_j -
(b^{ij}l_i)_{j}\big] - \Psi \nonumber \\
\ns&\ds   = \l^2 c_1^2 (2t - T)^2 + 4\l c_1 + \l c_0 -
\sum_{i,j}b^{ij}l_i l_j \\
\ns&\ds = \l^2  \Big[ c_1^2(2t - T)^2 -
\sum_{i,j}b^{ij}d_i
d_j \Big] + O(\l).
\end{array}
\end{equation}
By the definition of $B$, we see that
\begin{equation}\label{B1}
\begin{array}{ll}\ds
B \3n& \ds= A\Psi + (Al_t)_t -
\sum_{i,j}(Ab^{ij}l_i)_j + \frac{1}{2}\sum_{i,j}\big[ \Psi_{tt} -
(b^{ij}\Psi_i)_j\big]
 \\
\ns&\ds = 2Al_{tt} - \l c_0 A -\sum_{i,j}b^{ij}l_iA_j +
A_tl_t -\frac{1}{2}\sum_{i,j}\sum_{i',j'}\big[
b^{ij}(b^{i'j'}l_{i'})_{j'i}\big]_{j} \\
\ns&\ds  = 2\l^3 \Big[-2c_1^3(2t-T)^2 + 2c_1\sum_{i,j}b^{ij}d_id_j
\Big]
-\l^3c_0c_1^2(2t-T)^2 + \l^3c_0\sum_{i,j}b^{ij}d_id_j \\
\ns&\ds \q +
\l^3\sum_{i,j}\sum_{i',j'}b^{ij}d_i(b^{i'j'}d_{i'}d_{j'})_{j} -
4\l^3c_1^3(2t-T)^2 + O(\l^2)  \\
\ns&\ds  = (4c_1+c_0)\l^3 \sum_{i,j}b^{ij}d_id_j  +
\l^3\sum_{i,j}\sum_{i',j'}b^{ij}d_i(b^{i'j'}d_{i'}d_{j'})_{j} \\
\ns&\ds \q - (8c_1^3 + c_0c_1^2)\l^3(2t-T)^2 +
O(\l^2).
\end{array}
\end{equation}
Now we estimate
$\ds\sum_{i,j}\sum_{i',j'}b^{ij}d_i(b^{i'j'}d_{i'}d_{j'})_{j}
$. From Condition \ref{condition of d}, we get
that
\begin{equation}\label{bijdidj}
\begin{array}{ll}\ds
\mu_0\sum_{i,j}b^{ij}d_id_j \leq \sum_{i,j}\sum_{i',j'} \Big[
2b^{ij'}(b^{i'j}d_{i'})_{j'} - b_{j'}^{ij}b^{i'j'}d_{i'}
\Big]d_id_j   \\
\ns\ds\qq\qq\qq = \sum_{i,j}\sum_{i',j'}\Big(
2b^{ij'}b^{i'j}_{j'}d_{i'} + 2b^{ij'}b^{i'j}d_{i'j' } -
b_{j'}^{ij}b^{i'j'}d_{i'} \Big)d_id_j  \\
\ns\ds \qq\qq\qq =\sum_{i,j}\sum_{i',j'}\Big(
2b^{ij'}b^{i'j}_{j'}d_{i'}d_id_j + 2b^{ij'}b^{i'j}d_{i'j' }d_id_j
- b_{j'}^{ij}b^{i'j'}d_{i'}d_id_j \Big) \\
\ns\ds \qq\qq\qq =
\sum_{i,j}\sum_{i',j'}\Big(b^{i'j'}b_{j'}^{ij}d_{i'}d_id_j +
b^{ij}b^{i'j'}d_{i'j}d_id_j + b^{ij}b^{i'j'}d_{j'j}d_id_{i'} \Big)  \\
\ns\ds \qq\qq\qq =
\sum_{i,j}\sum_{i',j'}b^{ij}d_i(b^{i'j'}d_{i'}d_{j'})_j.
\end{array}
\end{equation}
From \eqref{B1} and \eqref{bijdidj}, by Condition \ref{condition2}, we
obtain that
\begin{eqnarray*}
 B \3n&\geq& \2n\ds\l^3 (4c_1+c_0)\sum_{i,j}b^{ij}d_id_j + \l^3 \mu_0
\sum_{i,j} b^{ij}d_id_j - (8c_1^3
+ 2c_0c_1^2)\l^3(2t-T)^2 + O(\l^2)  \\
 &\geq&\2n \ds  \l^3 (4c_1+c_0)\sum_{i,j}b^{ij}d_id_j + \l^3 \mu_0
\sum_{i,j} b^{ij}d_id_j - 2c_1^2(4c_1
+ c_0)\l^3T^2 + O(\l^2) \\
&\geq&\2n \ds 2(4c_1+c_0)\l^3 \(\sum_{i,j}b^{ij}d_id_j - c_1^2T^2\) + O(\l^2)  \\
& =&\2n\ds 2(4c_1+c_0)\l^3 (4R_1^2 - c_1^2T^2) + O(\l^2).
\end{eqnarray*}
Then we know that there exists a $\l_0 > 0$ such that for any $\l
\geq \l_0$, we have that
\begin{equation}\label{B ine}
Bv^2 \geq 8c_1  (4R_1^2 - c_1^2T^2) \l^3v^2.
\end{equation}

Since
$$
v(0,x)=\theta(0,x)y(0,x) = 0
$$
and
$$
v_t(0,x)=\theta_t(0,x)y(0,x)+\th(0,x)
y_t(0,x)=0,
$$
we know that at time
$t=0$, it holds that
$$
\sum_{i,j}b^{ij}l_t v_i v_j -
2\sum_{i,j}b^{ij}l_iv_jv_t + l_t v_t^2
- \Psi v_t v + \Big( Al_t +
\frac{\Psi_t}{2}\Big)v^2 =0.
$$
By a similar reason, we see that at time $t=T$,
$$
\sum_{i,j}b^{ij}l_t v_i v_j -
2\sum_{i,j}b^{ij}l_iv_jv_t + l_t v_t^2
- \Psi v_t v + \Big( Al_t +
\frac{\Psi_t}{2}\Big)v^2 = 0.
$$

{\bf Step 2.} Integrating
(\ref{hyperbolic2}) in $Q$, taking
expectation in $\O$ and by the argument
above, we obtain that
\begin{equation}\label{bhyperbolic31}
\begin{array}{ll}\ds
 \q\mathbb{E}\int_Q \theta\Big\{\Big( -2l_t v_t +
2\sum_{i,j}b^{ij}l_iv_j + \Psi v \Big)
\Big[ dy_t - \sum_{i,j}(b^{ij}y_i)_jdt \Big] - \theta l_t (dy_t)^2\Big\}dx  \\
\ns\ds\,\,\,\,\,\,\, + \l
\mathbb{E}\int_{\Si}\sum_{i,j}\sum_{i',j'}\big(
2b^{ij}b^{i'j'}d_{i'}v_i v_{j'} -
b^{ij}b^{i'j'}d_i v_{i'}v_{j'} \big)\nu^j d\Si  \\
\ns\ds \geq C \mathbb{E}\int_Q
\theta^2\Big[\big( \l  v_t^2 + \l
|\nabla v|^2  \big) +   \l^3 v^2 \]dxdt
+ \mathbb{E}\int_Q\Big(-2l_tv_t +
2\sum_{i,j}b^{ij}l_iv_j + \Psi v\Big)^2
dxdt.
\end{array}
\end{equation}
Since $y=0$ on $\Si$, $P$-a.s., from \eqref{def
gamma0}, we have
\begin{eqnarray}\label{bhyperbolic32}
&\,&\q\mathbb{E}\int_{\Si}\sum_{i,j}\sum_{i',j'}\Big(
2b^{ij}b^{i'j'}d_{i'}v_i v_{j'} - b^{ij}b^{i'j'}d_i v_{i'}v_{j'}
\Big)\nu^j d\Si  \nonumber\\
&\,&\ds = \mathbb{E}\int_{\Si}\sum_{i,j}\sum_{i',j'}\Big(
2b^{ij}b^{i'j'}d_{i'}\frac{\pa v}{\pa \nu}\nu^i \frac{\pa v}{\pa
\nu}\nu^{j'} - b^{ij}b^{i'j'}d_i \frac{\pa v}{\pa
\nu}\nu^{i'}\frac{\pa v}{\pa \nu}\nu^{j'}
\Big)\nu^j d\Si  \nonumber\\
&\,&\ds = \mathbb{E}\int_{\Si}\Big( \sum_{i,j}b^{ij}\nu^i \nu^j
\Big)\Big( \sum_{i',j'}b^{i'j'}d_{i'}\nu^{j'} \Big)\Big|\frac{\pa
v}{\pa \nu}\Big|^2d\Si\\
&\,&\ds =  \mathbb{E}\int_{\Si}\Big( \sum_{i,j}b^{ij}\nu^i \nu^j
\Big)\Big( \sum_{i',j'}b^{i'j'}d_{i'}\nu^{j'} \Big)\Big|\theta\frac{\pa
y}{\pa \nu} + y\frac{\pa\theta}{\pa\nu}\Big|^2d\Si \nonumber\\
&\,&\ds = \mathbb{E}\int_{\Si}\Big(
\sum_{i,j}b^{ij}\nu^i \nu^j \Big)\Big(
\sum_{i',j'}b^{i'j'}d_{i'}\nu^{j'}
\Big)\theta^2\Big|\frac{\pa y}{\pa
\nu}\Big|^2d\Si \nonumber\\
&\,&\ds \leq \mathbb{E}\int_{\Si_0}\Big(
\sum_{i,j}b^{ij}\nu^i \nu^j \Big)\Big(
\sum_{i',j'}b^{i'j'}d_{i'}\nu^{j'}
\Big)\theta^2\Big|\frac{\pa y}{\pa
\nu}\Big|^2d\Si.\nonumber
\end{eqnarray}

From \eqref{bhyperbolic31} and \eqref{bhyperbolic32}, we obtain that
\begin{equation}\label{bhyperbolic3}
\begin{array}{ll}
\ds\q\mathbb{E}\int_Q \theta\Big\{\Big( -2l_t v_t +
2\sum_{i,j}b^{ij}l_iv_j + \Psi v \Big)
\Big[ dy_t - \sum_{i,j}(b^{ij}y_i)_jdt \Big] - \theta l_t (dy_t)^2\Big\}dx  \\
\ns\ds\,\,\,\,\,\,\, + \l
\mathbb{E}\int_{\Si_0}\Big( \sum_{i,j}b^{ij}\nu^i
\nu^j \Big)\Big(
\sum_{i',j'}b^{i'j'}d_{i'}\nu^{j'}
\Big)\Big|\frac{\pa y}{\pa
\nu}\Big|^2d\Si  \\
\ns\ds \geq C \mathbb{E}\int_Q \Big[\theta^2\Big( \l  y_t^2 + \l
|\nabla y|^2  \Big) +   \l^3 \theta^2 y^2 \]dxdt + \mathbb{E}\int_Q\Big(-2l_tv_t +
2\sum_{i,j}b^{ij}l_iv_j + \Psi
v\Big)^2 dxdt.
\end{array}
\end{equation}

Since  $y$ solves the equation \eqref{system2},  we know that
\begin{eqnarray}\label{bhyperbolic4}
&\,&\q\mathbb{E}\int_Q \theta\Big\{\Big(
-2l_t v_t + 2\sum_{i,j}b^{ij}l_iv_j +
\Psi v \Big) \Big[ dy_t -
\sum_{i,j}(b^{ij}y_i)_jdt \Big] - \theta l_t (dy_t)^2\Big\}dx   \nonumber\\
&\,&\ds= \mathbb{E}\int_Q
\theta\Big\{\Big( -2l_t v_t +
2\sum_{i,j}b^{ij}l_iv_j + \Psi v \Big)
\big[ b_1 y_t +  b_2\cd\nabla y + b_3 y
+ \chi  f + \a\big]
\\ &\,&\ds\qq\qq - \theta l_t (b_4 y + \chi
g)^2\Big\}dxdt  \nonumber
\\
&\,&\ds\leq  C\Big\{ \mathbb{E}\int_Q
\theta^2 \Big[ b_1y_t + b_2\cd\nabla y
+ b_3 y + \chi f + \a \Big]^2 + \l
\theta^2(b_4 y + \chi  g)^2 \Big\}dxdt \nonumber
\\
&\,&\ds \qq +\mathbb{E}\int_Q\(-2l_t v_t
+ \sum_{i,j}b^{ij}l_i v_j +
\Psi v \)^2dxdt  \nonumber\\
&\,&\ds \leq C\bigg\{\mathbb{E}\int_Q
\theta^2(f^2 + \a^2 + \l g^2)dxdt +
|b_1|^2_{L^{\infty}_{\cF}(0,T;L^{\infty}(G))}
\mathbb{E}\int_Q \theta^2 y_t^2 dxdt +
\mathbb{E}\int_Q \theta^2 b_3^2
y^2 dxdt  \nonumber
\\
&\,&\ds \qq\;\, +
|b_2|^2_{L^{\infty}_{\cF}(0,T;L^{\infty}(G,\mathbb{R}^n))}  \mathbb{E}\int_Q\theta^2
|\nabla y|^2dxdt +  \l
|b_4|_{L^{\infty}_{\cF}(0,T;L^{\infty}(G))}^2\mathbb{E}\int_Q
\theta^2 y^2 dxdt   \bigg\}\nonumber
 \\
&\,&\ds \q +\mathbb{E}\int_Q\(-2l_t v_t +
\sum_{i,j}b^{ij}l_i v_j + \Psi v
\)^2dxdt.\nonumber
\end{eqnarray}

Recalling the definition of $r_2$ in \eqref{r1r2}, and using successively H\"{o}lder's and Sobolev's inequalities, we get
\begin{equation}\label{bhyperbolic4.1}
|b_3\theta y|^2_{L^2_\cF(0,T;L^2(G))} \leq r_2 |\theta y|^2_{L^2_\cF(0,T;L^{s}(G))} \leq r_2|\theta y|^2_{L^2_\cF(0,T;H^{n/p})} \mbox{ for } \frac{1}{p} + \frac{1}{s} = \frac{1}{2}.
\end{equation}
For any $F\in L^2(\O,\cF_T,P;H^1(\dbR^n))$, by H\"{o}lder's inequality, one has
$$
\begin{array}{ll}\ds
|F|^2_{L^2(\O,\cF_T,P;H^{n/p}(\dbR^n))} = \mE\int_{\dbR^n}(1+|\xi|^2)^{n/p}|\hat F(\xi)|^{2n/p}|\hat F(\xi)|^{2(1-n/p)}d\xi \\
\ns\ds \hspace{3.8cm} \leq |F|^{2n/p}_{L^2(\O,\cF_T,P;H^1(\dbR^n))}|F|^{2(1-n/p)}_{L^2(\O,\cF_T,P;L^2(\dbR^n))}.
\end{array}
$$
Hence, we know that there is a constant $C>0$ such that for any $\widetilde F \in L^2(\O,\cF_T,P;H_0^1(G))$, we have
$$
|\widetilde F|^2_{L^2(\O,\cF_T,P;H^{n/p}(G))} \leq C|\widetilde F|^{2n/p}_{L^2(\O,\cF_T,P;H_0^1(G))}|\widetilde F|^{2(1-n/p)}_{L^2(\O,\cF_T,P;L^2(G))}.
$$
Therefore, there is a constant $C>0$ such that for any $\overline F \in L^2_\cF(0,T;H_0^1(G))$, it holds that
$$
|\overline F|^2_{L^2_\cF(0,T;H^{n/p}(G))} \leq C|\overline F|^{2n/p}_{L^2_\cF(0,T;H_0^1(G))}|\overline F|^{2(1-n/p)}_{L^2_\cF(0,T;L^2(G))}.
$$
This, together with the inequality \eqref{bhyperbolic4.1}, implies that
\begin{equation}\label{bhyperbolic4.2}
\begin{array}{ll}\ds
|b_3\theta y|^2_{L^2_\cF(0,T;L^2(G))}\leq C|b_3\theta y|^{2n/p}_{L^2_\cF(0,T;H_0^1(G))}|b_3\theta y|^{2(1-n/p)}_{L^2_\cF(0,T;L^2(G))}\\
\ns\ds \hspace{3.15cm} \leq \e\l|b_3\theta y|^{2}_{L^2_\cF(0,T;H_0^1(G))} + C(\e)r_2^{ 2p/(p-n)}\l^{-n/(p-n)}|b_3\theta y|^{2}_{L^2_\cF(0,T;L^2(G))},
\end{array}
\end{equation}
where $\e$ is small enough and $C(\e)$ depends on  $\e$.

Taking $\l_2 = C(r_1^2 +
r_2^{\frac{1}{3/2-n/p}}+1)\geq\max\{\l_0,\l_1\}$,
combining \eqref{bhyperbolic3},
\eqref{bhyperbolic4} and \eqref{bhyperbolic4.2},
for any $\l \geq \l_2$, we have that
\begin{equation}\label{bhyperbolic5}
\begin{array}{ll}
\ds\q  C \l \mathbb{E}\int_{\Si_0}\theta^2\Big(
\sum_{i,j}b^{ij}\nu_i \nu_j \Big)\Big(
\sum_{i',j'}b^{i'j'}d_{i'}\nu_{j'}
\Big)\Big|\frac{\pa
y}{\pa \nu}\Big|^2d\Si  + C\mathbb{E}\int_Q \theta^2(f^2 + \a^2 + \l g^2)dxdt \\
\ns\ds \geq \mathbb{E}\int_Q \theta^2 \(  \l  y_t^2 + \l
|\nabla y|^2   +   \l^3 y^2  \)dxdt.
\end{array}
\end{equation}
Recalling the property of $\chi$(see \eqref{chi}) and $y=\chi z$, from \eqref{bhyperbolic5}, we find
\begin{equation}\label{bhyperbolic6}
\begin{array}{ll}
\ds\q  C \l \mathbb{E}\int_{\Si_0}\Big(
\sum_{i,j}b^{ij}\nu_i \nu_j \Big)\Big(
\sum_{i',j'}b^{i'j'}d_{i'}\nu_{j'}
\Big)\theta^2\Big|\frac{\pa
z}{\pa \nu}\Big|^2d\Si  + C\mathbb{E}\int_Q \theta^2(f^2  + \l g^2)dxdt \\
\ns\ds \q + C(r_1 + 1)\[\mathbb{E}\int_{0}^{\frac{T}{2}-\e_1T}\!\!\!\!\int_G \theta^2 (z_t^2 \!+\! |\nabla z|^2 \!+\! z^2)dxdt \!+\! \mathbb{E}\int_{\frac{T}{2}+\e_1T}^{ T}\int_G \theta^2 (z_t^2 + |\nabla z|^2 + z^2)dxdt\] \\
\ns\ds \geq \mathbb{E}\int_{\frac{T}{2}-\e_0
T}^{\frac{T}{2}+\e_0 T}\int_G \theta^2 \(  \l
z_t^2 + \l |\nabla z|^2   +   \l^3  z^2  \)dxdt.
\end{array}
\end{equation}
Combining \eqref{en esti} and \eqref{bhyperbolic6}, we know that there is a $\l_3=C(r_1^2 + r_2^{\frac{1}{3/2-n/p}}+1)\geq \l_2$ such that for all $\l\geq \l_3$, it holds that
\begin{equation}\label{bhyperbolic7}
\begin{array}{ll}
\ds\q  C \l \mathbb{E}\int_{\Si_0}\theta^2\Big(
\sum_{i,j}b^{ij}\nu_i \nu_j \Big)\Big(
\sum_{i',j'}b^{i'j'}d_{i'}\nu_{j'}
\Big)\Big|\frac{\pa
z}{\pa \nu}\Big|^2d\Si  + C\mathbb{E}\int_Q \theta^2(f^2  + \l g^2)dxdt  \\
\ns\ds \geq e^{-\frac{T^2}{4}\l_3}\mathbb{E} \int_G  \( z_1^2 +
|\nabla z_0|^2  \)dxdt.
\end{array}
\end{equation}
Taking $\l=\l_3$, we obtain that
\begin{equation}\label{bhyperbolic8}
\begin{array}{ll}
\ds\q  C e^{\l_3 R_1^2}\[\mathbb{E}\int_{\Si_0}
\Big|\frac{\pa
z}{\pa \nu}\Big|^2d\Si  +  \mathbb{E}\int_Q  (f^2  + g^2)dxdt \] \\
\ns\ds \geq e^{-\frac{T^2}{4}\l_3}\mathbb{E} \int_G  \big( z_1^2 +
|\nabla z_0|^2   \big)dxdt.
\end{array}
\end{equation}
This leads to the inequality \eqref{obser esti2} immediately.

\section{Proof of Theorem \ref{inobser}}

\q\,\, This section is devoted to a
proof of  Theorem \ref{inobser}.

\vspace{0.2cm}

{\it Proof of Theorem \ref{inobser}}\,:
Let $h_0\in C^1(\overline G;
\mathbb{R}^n)$ such that $h_0 = \nu$ on
$\G$, and let $\rho\in C^2(\overline
G;[0,1])$ such that

\vspace{0.1cm}

\begin{equation}\label{rho}
\left\{
\begin{array}{ll}
\ds \rho =1 &\mbox{ in }  \cO_{\frac{\d}{3}}(\G_0),\\
\ns\ds \rho = 0 &\mbox{ in } G\setminus \cO_{\frac{\d}{2}}(\G_0).
\end{array}
\right.
\end{equation}

Let $h = \rho\theta^2 h_0$ in the equality
\eqref{equality hidden1}, noting that
$y_j = \frac{\pa y}{\pa\nu}\nu^j$ on
$\Si$, by integrating by parts, we see
\begin{equation}\label{sta1}
\begin{array}{ll}
\ds \q\mE\int_{\Si} \(\sum_{i,j=1}^n b^{ij}\nu^i\nu^j \)\rho\theta^2 \Big|\frac{\pa y}{\pa \nu}\Big|^2 d\G dt\\
\ns\ds = \mE\int_Q \sum_{i=1}^n\Big[
2(h\cd\nabla y)\sum_{j=1}^n b^{ij}y_{j}
+ h^i\Big(y_t^2 - \sum_{j,k=1}^n
b^{jk}y_{j} y_{k} \Big) \Big]_{i}dxdt\\
\ns\ds = -\mE\int_Q 2 \Big[\Big(dy_t -
\sum_{i,j=1}^n
(b^{ij}y_{i})_{j}dt\Big) h \cd
\nabla y - d(y_t h\cd \nabla
y)  + y_t h_t \n ydt\\
\ns\ds \qq\qq\q - \sum_{i,j,k=1}^n
b^{ij}y_{i} y_{k} h^k_{j}dt  -
\div (h) y_t^2dt + \sum_{i,j=1}^n
y_{j} y_{i} \div(b^{ij}h)\Big] dxdt\\
\ns\ds \leq
C\Big\{\frac{1}{\l}\mE\int_Q\theta^2 \big(
b_1y_t +  b_2\cd\nabla y + b_3 y + \chi
f + \a \big)^2 dxdt \\
\ns\ds \qq \; + \l\mE\int_0^T\int_{
\cO_{\frac{\d}{2}}(\G_0)}
\theta^2(y_t^2 + |\nabla
y|^2)dxdt\Big\}.
\end{array}
\end{equation}

Now let us deal with the term
$\ds\mE\int_0^T\int_{ \cO_{\frac{\d}{2}}(\G_0)}
 \theta^2 |y_t|^2 dxdt$. Let $\rho_1 \in C^2(\overline G;[0,1])$ satisfying that
 \vspace{0.1cm}
$$
\left\{
\begin{array}{ll}
\ds \rho_1 = 1 &\mbox{ in }  \cO_{\frac{\d}{2}}(\G_0),\\
\ns\ds \rho_1 = 0 &\mbox{ in } G\setminus \cO_{\d}(\G_0).
\end{array}
\right.
$$
Put $\eta = \rho_1^2 \theta^2$. By
virtue of that $y$ solves the equation
\eqref{system2}, we have
\begin{eqnarray}\label{sta2}
&\,&
\q\mE\int_Q \eta y \big(b_1 y_t + b_2\cd\nabla y + b_3 y + f + \a\big)dxdt \nonumber\\
&\,&\ds = \mE\int_Q \eta y \[dy_t - \sum_{i,j=1}^n (b^{ij}y_{i})_j dt\] dx\\
&\,&\ds   = -\mE\int_Q \[ y_t (\eta_t y + \eta y_t) \]dxdt + \mE\int_Q \eta \sum_{i,j=1}^n b^{ij}y_i y_j dxdt \nonumber\\
&\,&\ds  \q + \mE\int_Q  y\sum_{i,j=1}^n
b^{ij}y_i \eta_j dxdt,\nonumber
\end{eqnarray}
 this implies that
\begin{equation}\label{sta3}
\begin{array}{ll}
\ds\q \mE\int_0^T\int_{\cO_{\frac{\d}{2}}(\G_0)} \theta^2 |y_t|^2 dxdt\\
\ns\ds \leq C\Big\{ \frac{1}{\l^2}\mE \int_Q
\theta^2 \big[ b_1y_t +  b_2\cd\nabla y
+ b_3 y + \chi f + \a \big]^2 dxdt \\
\ns\ds \qq +
\mE\int_0^T \int_{\cO_{\d}(\G_0)} \theta^2 (\l^2
y^2 + |\nabla y|^2)dxdt \Big\}.
\end{array}
\end{equation}
From \eqref{bhyperbolic5}, \eqref{sta2}
and \eqref{sta3}, we get that there is
a $\l_4=C(r_1^2 +
r_2^{\frac{1}{3/2-n/p}}+1)>0$   such that for any $\l \geq
\max\{\l_3, \l_4\}$, it holds that
\begin{equation}\label{sta4}
\begin{array}{ll}\ds
\q e^{-\l\frac{T^2}{4}}\mathbb{E} \int_G
(|\nabla z_0|^2  +
|z_1|^2 ) dx
\\ \ns\ds\leq C\mE\int_0^T\int_{\cO_{\d}(\G_0)} \theta^2 (\l^3 z^2 + \l |\n z|^2)dxdt.
\end{array}
\end{equation}

Since
$$
e^{-\l\frac{T^2}{4}} \leq \theta(t,x)
\leq e^{\l\frac{T^2}{4}},
$$
we see that
\begin{equation}\label{sta5}
\mathbb{E} \int_G  \big(|\nabla z_0|^2
+ |z_1|^2\big)dx   \leq C
e^{\frac{\l T^2}{2}}
\mE\int_0^T\int_{\cO_{\d}(\G_0)} (|\n z|^2+z^2)
dxdt.
\end{equation}
This, together with Poincar\'{e}'s
inequality, implies the inequality
\eqref{inobser esti2} immediately.

\section{A state observation problem }\label{Sec app}

This section is addressed to a state observation problem for semilinear stochastic hyperbolic equations.
Let
$$
F(\eta,\varrho,\zeta):
\mathbb{R}^1\t\mathbb{R}^1\t\mathbb{R}^n\to\
\mathbb{R}^1
$$
and
$$
K(\eta):
\mathbb{R}^1\to\ \mathbb{R}^1
$$
be two known nonlinear functions.
Consider the following semilinear stochastic
hyperbolic equation
\begin{equation}\label{4.12-eq1}
\left\{
\begin{array}{ll}
\ds dw_t-
\sum_{i,j}(b^{ij}w_i)_jdt=F(w,w_t,\n
w)dt + K(w)
dB(t)&\mbox{ in }Q,\\
\ns\ds w=0&\mbox{ on }\Si,\\
\ns\ds w(0)=w_0,\q w_t(0)=w_1&\mbox{ in
}G,
\end{array}
\right.
\end{equation}
where the initial data $(w_0,w_1)\in
L^2(\O,\cF_0,P;$ $H_0^1(G)\times
L^2(G))$ are unknown random variables.

We put the following assumption:

{\bf (AS)} The nonlinear functions
$F(\cd,\cd,\cd)$ and $K(\cd)$ satisfy
the following:
\begin{enumerate}
  \item $$
\ba{l@{\,}l}
\ds|F(\eta_1,\varrho,\zeta)-F(\eta_2,\varrho,\zeta)|&\le
L(1+|\eta_1|^{p-1}+|\eta_2|^{p-1})|\eta_1-\eta_2|
\\ \ns &\ds\forall\ \eta_1, \eta_2,\varrho\in
\mathbb{R}^1,\ \zeta\in \mathbb{R}^n
 \ea
$$
 with $1\le p\le \frac{n}{n-2}$ if
$n\ge3$; $1\le p<\infty$ if $n=1, 2$, for some constant $L>0$;
\item 
$$
\begin{array}{ll}
\ds |F(\eta,\varrho_1,\zeta_1)-F(\eta,\varrho_2,\zeta_2)|\le
L(|\varrho_1-\varrho_2|+|\z_1-\z_2|) \\
\ns \ds\qq\qq\qq\forall\
(\eta,\varrho_i,\z_i)\in
\mathbb{R}^1\t\mathbb{R}^1\t\mathbb{R}^n,\
i=1,2, \\
\ns\ds |F(0,\varrho,\z)|\le
L(|\varrho|+|\z|)\ \ \forall \ (v,\z)\in
\mathbb{R}^1\t\mathbb{R}^n,\\
\ns\ds  |K(\eta_1)-K(\eta_2)| \leq L|\eta_1-\eta_2|\forall\ \eta_1, \eta_2, \in
\mathbb{R}^1
\end{array}
$$
for some constant $L>0$;
  \item for any given initial data
$(w_0,w_1)\in
L^2(\O,\cF_0,P;H_0^1(G)\times L^2(G))$,
\eqref{4.12-eq1} admits a unique
solution $w=w(\cd\, ;w_0,w_1)\in H_T$
(the solution of \eqref{4.12-eq1}  is
defined similarly to
the one of \eqref{system1}).
\end{enumerate}

Here since we do not introduce any sign
condition on the nonlinear functions
$F(\cd,\cd,\cd)$ and $K(\cd)$, the global
existence of a solution to
\eqref{4.12-eq1} is not guaranteed.
This is why we need to impose the third
assumption in {\bf (AS)}.

The state observation problem
associated to the equation
\eqref{4.12-eq1} is as follows.

\vspace{0.1cm}
\begin{itemize}

\item {\bf Identifiability}. Is the solution
$w\in H_T$ (to \eqref{4.12-eq1})
determined uniquely by the observation
$\ds\frac{\pa
w}{\pa\nu}\Big|_{(0,T)\times
\G_0}$(\resp $ w|_{(0,T)\times
\cO_\d(\G_0)}$)?

\vspace{0.1cm}

\item {\bf Stability}. Assume that two
solutions $w$ and $\hat w$ (to the equation
\eqref{4.12-eq1}) are given. Let
$\ds\frac{\pa
w}{\pa\nu}\Big|_{(0,T)\times \G_0}$(\resp $w|_{(0,T)\times\cO_\d(\G_0)}$) and
$\ds\frac{\pa \hat
w}{\pa\nu}\Big|_{(0,T)\times \G_0}$(\resp $\hat w|_{(0,T)\times\cO_\d(\G_0)}$) be
the corresponding observations. Can we
find a positive constant $C$ such that
$$
| w-\hat w | \leq C\|\!\| \frac{\pa
w}{\pa\nu}-\frac{\pa \hat w}{\pa\nu}
\|\!\|\(\resp | w-\hat w | \leq C|\!|
 w- \hat w |\!|\),
$$
with appropriate norms in both sides?

\vspace{0.1cm}

\item {\bf Reconstruction}. Is it possible to
reconstruct $w\in H_T$ to
\eqref{4.12-eq1}, in some sense,  from
the observation $\ds\frac{\pa
w}{\pa\nu}\Big|_{(0,T)\times
\G_0}$(\resp $ w|_{(0,T)\times
\cO_\d(\G_0)}$)?

\end{itemize}

The state observation problem for
systems governed by deterministic
partial differential equations is
studied extensively (see
\cite{Kli,Li1,Yamamoto} and the rich
references therein). However, the
stochastic case attracts very little
attention. To our best knowledge,
\cite{Zhangxu3,Luqi4} are the only two
published papers addressing this topic.
In \cite{Luqi4}, the author studied the
state observation problem for
stochastic Schr\"{o}dinger equations
via the Carleman estimate for the
equation. In \cite{Zhangxu3}, the
author addressed the state observation
problem for stochastic wave equations
and  proved the following result
\begin{equation}\label{4.12-eq2}
\begin{array}{ll}\ds
|(w(t)-\hat w(t), w_t(t)-\hat
w_t(t))|_{L^2(\O,\cF_t,P;H_0^1(G)\t
L^2(G))} \\
\ns
 \ds\le e^{Ct^{-1}}\wt C\Big|\frac{\pa
w}{\pa\nu}-\frac{\pa \hat
w}{\pa\nu}\Big|_{L^2_{\cF}(0,T;L^2(\G_0))}
\q \mbox{ for any } t>0.
\end{array}
\end{equation}
Obviously, one cannot let $t=0$ in
\eqref{4.12-eq2}, which means that the initial state cannot be obtained from the observation. In this paper, by
means of Theorem \ref{observability},
we can give positive answers to the
above first and second questions, that
is, we prove that the whole state can
be observed by the boundary or internal
observation.

First,  thanks to the Sobolev embedding
theorem and the conditions on $F(\cd,\cd,\cd)$ and $K(\cd)$, we know
$$
F(w,w_t,\n w)\in
L^2_{\cF}(0,T;L^2(G)),\quad  K(w)\in
L^2_{\cF}(0,T;L^2(G))
$$
for any $w\in H_T$. Thus, by
Proposition~\ref{hidden r}, we know
$\frac{\pa w}{\pa\nu}\in
L^2_\cF(0,T;L^2(G_0))$. Now, we define
two nonlinear maps $\cM_1$ and $\cM_2$ as follows:
$$
\left\{
\begin{array}{ll}\ds \cM_1:\
L^2(\O,\cF_0,P;H_0^1(G)\times
L^2(G))\to
L^2_{\cF}(0,T;L^2(\G_0)),\\
\ns\ds \cM_1(w_0,w_1)= {\frac{\pa
w}{\pa \nu}}\Big|_{(0,T)\times
\G_0},
\end{array}
\right.
$$
$$
\left\{
\begin{array}{ll}\ds \cM_2:\
L^2(\O,\cF_0,P;H_0^1(G)\times
L^2(G))\to
L^2_{\cF}(0,T;L^2(\cO_\d(\G_0))),\\
\ns\ds \cM_2(w_0,w_1)= {\n
w}\big|_{(0,T)\times \cO_\d(\G_0)},
\end{array}
\right.
$$
where $w$ solves the equation
\eqref{4.12-eq1}.

\vspace{0.3cm}

We have the following result.

\begin{theorem}\label{th2}
Let   Condition (\ref{condition of d})
and Condition (\ref{condition2}) be
satisfied. There exists a constant $\wt C=\wt
C(L,T,G,(b^{ij})_{1\leq i,j\leq n},\G_0,\d)>0$ such that for any initial data
$(w_0,w_1), (\hat w_0, \hat w_1)\in
L^2(\O,\cF_0,P;H_0^1(G)\times L^2(G))$,
it holds that
\begin{equation}\label{th2eq1}
|(w_0-\hat w_0, w_1-\hat
w_2)|_{L^2(\O,\cF_0,P;H_0^1(G)\times
L^2(G))} \le \wt C|\cM_1(w_0,w_1)
-\cM_1(\hat w_0,\hat
w_1)|_{L^2_{\cF}(0,T;L^2(\G_0))}
\end{equation}
and
\begin{equation}\label{th2eq2}
|(w_0-\hat w_0, w_1-\hat
w_2)|_{L^2(\O,\cF_0,P;H_0^1(G)\times
L^2(G))} \le \wt C|\cM_2(w_0,w_1)
-\cM_2(\hat w_0,\hat
w_1)|_{L^2_{\cF}(0,T;L^2(\cO_\d(\G_0)))},
\end{equation}
where $\hat w=\hat w(\cd\, ;\hat
w_0,\hat w_1)\in H_{T}$ is the solution
to \eqref{4.12-eq1} with $(w_0,w_1)$
replaced by $(\hat w_0,\hat w_1)$.
\end{theorem}

\begin{remark}    Theorem~\ref{th2}
indicates that the state $w(t)$ of
\eqref{4.12-eq1} (for $t\in [0,T]$) can
be uniquely determined from the
observed boundary data $\ds{\frac{\pa
w}{\pa \nu}} \Big|_{(0,T)\t\G_0}$ or $\n w|_{(0,T)\t \cO_\d(\G_0)}$,
$P$-a.s., and continuously depends on
it. Therefore, we answer the first and
second questions for the state observation
problem of the system \eqref{4.12-eq1}
positively.
\end{remark}

{\it Proof of Theorem}~\ref{th2}: Set
$$
y=\hat w-w.
$$
It is easy to see that $y$ is a
solution of \eqref{4.12-eq1} with
$$
\left\{
\begin{array}{ll}\ds
b_1=\int_0^1\pa_\varrho(\hat
w,w_t+s(\hat w_t-w_t),\n w)ds, \q b_2=\int_0^1\pa_\zeta F(\hat w,\hat
w_t,\n w+s(\n \hat w-\n w))ds, \cr \ns
\ds b_3=\int_0^1\pa_\eta F(w+s(\hat w-w),
w_t ,\n w)ds,\q
b_4=\int_0^1\pa_\eta K(w+s(\hat
w-w))ds.
\end{array}\right.
$$
Then, the inequality \eqref{th2eq1}
follows from Theorem \ref{observability} and
the inequality \eqref{th2eq2} comes
from Theorem \ref{inobser}.
\endpf

As a direct consequence of Theorem
\ref{th2}, we have the following unique
continuation property for the equation
\eqref{system1}.

\begin{corollary}\label{thucp}
Let   Condition (\ref{condition of d})
and Condition (\ref{condition2}) be
satisfied. Assume that $f=g=0$ in $Q$,
$P$-a.s. If  a solution of the equation
\eqref{system1} satisfies $y=0$ in
$(0,T)\t O_\d(\G_0)$, $P$-a.s., then we
have that $y=0$ in $Q$, $P$-a.s.
\end{corollary}

\begin{remark}
The analogous result of Corollary
\ref{thucp} for deterministic
hyperbolic equations with nonsmooth
lower order terms was first obtained in
\cite{Zhangxu1}.
\end{remark}

Due to the need from Control/Inverse
Problems of partial differential
equations, the study of the global
unique continuation for partial
differential equations is very
active(see
\cite{Castro-Zuazua,Zhangxu1,Zhang-Zuazua}
and the references therein) in recent
years. Compared with the plentiful
studying of the unique continuation
property for partial differential
equations, the study  for stochastic
partial differential equations is cold
and cheerless. To the best of our
knowledge,
\cite{Zhangxu4,Zhangxu5,Luqi4,Luqi5} are the
only published articles which concern
this topic, and the above unique
continuation property for stochastic
hyperbolic equations has not been
presented in the literature.

Next, we consider the reconstruction of the state $w$. Denote by $\f$ the observation on $(0,T)\times \G_0$ and by
$$
W\=\{w\in H_T:\,w \mbox{ solves \eqref{4.12-eq1} for some initial data } (w_0,w_1)\in L^2(\O,\cF_0,P;H_0^1(G)\times
L^2(G))\}.
$$  Put
$$
J_1(w) = \mE\int_0^T\int_{\G_0} \Big| \frac{\pa w}{\pa\nu} - \f \Big|^2 d\G dt \q\mbox{ for } w\in W.
$$
Let $\tilde w$ be the state of \eqref{4.12-eq1} corresponding to the observation $\f$. Then, it is clear that the state $\tilde w$ satisfies that
$$
J_1(\tilde w) = \min_{w\in W} J_1(w) = 0.
$$
Hence, the construction of $\tilde w$ can be converted to the study of the following optimization problem
$$
{\rm (P_1)} \mbox{ Find a  } w\in W \mbox{ which minimizes } J_1(\cd).
$$

By a similar argument, we can show the construction of $\tilde w$ can be deduced to  the following optimization problem
$$
{\rm (P_2)} \mbox{ Find a  } w\in W \mbox{ which minimize } J_2(w)=\int_0^T\int_{G_0}|\n w - \psi|^2 dxdt,
$$
where $\psi$ is the internal observation.

To give  efficient algorithms to solve problem ($\rm P_1$) and ($\rm P_2$) is beyond the scope of this paper and will be studied in our forthcoming paper.

\section{Further comments and open problems}

There are plenty of open problems in
the topic of this paper. Some of them
are particularly relevant and could
need  new ideas and further
developments.

\begin{itemize}

\item {\bf Efficiency algorithm for the construction of the solution $w$ from the observation}

In this paper, we only answer the first
and the second questions in the state
observation problem. The third one is
still open. Due to the stochastic
feature, some efficient approaches for
hyperbolic equation(see \cite{Li1} for
example), become invalid. In the end of Section \ref{Sec app}, we show that it can be solved by studying an optimization problem.
In this context, it seems that
one may utilize the great many sharp
methods   from optimization theory to
study the construction of $(w_0,w_1)$.
However, thanks to the stochastic
setting, it seems that one cannot
simply mimic these methods. A detailed
study of this  interesting but
difficult problem is beyond the scope
of this paper.

\item{\bf Observability estimate and unique continuation property with less restrictive conditions}

In this paper, we prove the inequality
\eqref{obser esti2} and \eqref{inobser
esti2} under the Condition
\eqref{condition of d} and Condition
\eqref{condition2}.   It is well known
that a sharp sufficient condition for
establishing observability estimate for
deterministic hyperbolic equations with
time invariant lower order terms   is
that the triple $(G, \G_0,
T)$($(G,\cO_\d(\G_0),T)$) satisfies the
geometric optic condition introduced in
\cite{Bardos-Lebeau-Rauch1}. It would
be quite interesting and challenging to
extend this result to the stochastic
setting. However, there are
lots of things should be done before
solving this problem. For instance, the
propagation of singularities for
stochastic partial differential
equations, at least, for stochastic
hyperbolic equations, should be
established.

As we have pointed out in Remark
\ref{rmI}, it is more interesting to
get the following inequality
\begin{equation}\label{4.12-eq3}
\begin{array}{ll}\ds
\q|(z_0,z_1)|_{L^2(\O,{\cal F}_0, P;
L^2(G)\t H^{-1}(G))}
\\ \ns\ds \leq e^{C(r_1^2 + r_2^{\frac{1}{ 3/2 - n/p}}+1)} \Big(| z  |_{L^2_{\cal
 F}(0,T;L^2(\cO_\d(\G_0)))} + |f|_{L^2_{\cal
 F}(0,T;L^2(G))} + |g|_{L^2_{\cal
 F}(0,T;L^2(G))}\Big).
\end{array}
\end{equation}

For deterministic hyperbolic equations,
the inequality \eqref{4.12-eq3} can be
obtained by combing the global Carleman
estimate and the multiplier method(see
\cite{Fu-Yong-Zhang1} for example). If
one follows the method to study the
stochastic problem, one will meet some
real difficulty. In fact, as the
deterministic settings, from Theorem
\ref{observability}, by means of a
suitable choice of multiplier, one can
get
\begin{equation}\label{rem eq1}
\begin{array}{ll}
 \ds\q |(z_0,z_1)|_{L^2(\O,{\cal F}_0, P; H_0^1(G)\t
L^2(G))}\\
 \ns\ds \leq Ce^{C(r_1^2 + r_2^{\frac{1}{ 3/2 - n/p}}+1)} \Big(\mE\int_0^T\int_{\cO_\d(\G_0)}\big(z_t^2 + z^2\big) dxdt + |f|_{L^2_{\cal
 F}(0,T;L^2(G))} + |g|_{L^2_{\cal
 F}(0,T;L^2(G))}\Big).
\end{array}
\end{equation}
Then, by employing the
Compact/Uniqueness argument, we can
eliminate the term ``$ z^2$" in the
right hand side of \eqref{rem eq1} for
deterministic case. A key point in the
Compact/Uniqueness argument is the fact
that $H^1((0,T)\times G)$ is compactly
imbedded into $L^2((0,T)\times G)$.
However,  the corresponding result is
not true in the stochastic settings.
One can easily show that even
$H_T$(recall that $H_T$ is given as in
\eqref{HT}) is not compact embedded in
$L^2_\cF(0,T;L^2(G))$.   The missing of
compactness leads to new difficulty for
establishing internal observability
estimate for stochastic hyperbolic
equations.

Under the Condition \eqref{condition of
d} and Condition \eqref{condition2},
$y=0$ in $Q$, $P$-a.s., provided that
$y=0$ in $(0,T)\times \cO_\d(\G_0)$.  Compared
to the classical unique continuation
result for deterministic hyperbolic
equations (see \cite{Tataru} for
example), the conditions in this paper
is very restrictive. It would be quite
interesting but maybe challenging to
prove whether these results in
\cite{Tataru}  is true or not for
stochastic hyperbolic equations.

\item {\bf Some other inverse problems for stochastic hyperbolic equations}

    In this paper, we show that the state can be uniquely determined by the observation via Carleman estimate. For deterministic partial differential equations, there are lots of other interesting inverse problems solved by some methods based on Carleman estimate. For example,  the multidimensional coefficient/source inverse problems with single measurement data. Both the global uniqueness and global stability are obtained by some methods which are mainly based on Carleman estimate. There are so many works in this topic. Hence, we do  not list them comprehensively and we refer the readers to two nice surveys \cite{Kli2,Yamamoto} and the rich references therein. One will meet substantially new difficulties in the study of inverse problems for stochastic partial differential equations. For instance, unlike the deterministic partial differential equations, the solution of a stochastic partial differential equations is usually non-differentiable with respect to the variable with noise (say, the time variable considered in this paper). Also, the usual compactness embedding result does not remain true for the solution spaces related to stochastic partial differential equations. Due to these new phenomenons, most of the powerful methods in \cite{Kli2,Yamamoto} cannot be applied to stochastic partial differential equations directly. In \cite{Luqi3}, the author studied an inverse source problem for stochastic parabolic equations involved in some special domain. It seems that it is hard to generalize the method in \cite{Luqi3} for the study of stochastic partial differential equations in general domains.

\end{itemize}




\end{document}